\documentclass[12pt]{amsart}
\usepackage[latin1]{inputenc}
\usepackage{amsxtra,amssymb,amsthm,amsmath,amscd,mathrsfs}
\usepackage[all]{xy}
\def \R {{\mathbb R}}

\def \Z {{\mathbb Z}}

\def\re{{\Re e\,}}

\def\le{\leqslant}
\def\ge{\geqslant}

\theoremstyle{plain}
\newtheorem{theorem}{Theorem}

\newtheorem{lemma}{Lemma}[section]
\newtheorem{corollary}{Corollary}

\theoremstyle{remark}

\theoremstyle{definition}

\numberwithin{equation}{section}

\begin{document}

\title[Power sums of Hecke eigenvalues and application]
{Power sums of Hecke eigenvalues and application} 
\author{J. Wu}
\address{Institut Elie Cartan Nancy (IECN)
\\
Nancy-Universit\'e, CNRS, INRIA
\\
Boulevard des Aiguillettes, B.P. 239
\\
54506 Van\-d\oe uvre-l\`es-Nancy
\\
France}
\email{wujie@iecn.u-nancy.fr}

\date{\today}

\begin{abstract}
We sharpen some estimates of Rankin on power sums of Hecke eigenvalues,
by using Kim \& Shahidi's recent results on higher order symmetric powers.
As an application, 
we improve Kohnen, Lau \& Shparlinski's lower bound 
for the number of Hecke eigenvalues of same signs.
\end{abstract}
\maketitle

\addtocounter{footnote}{1}

\section{Introduction}

\smallskip

Let $k\ge 2$ be an even ineteger and $N\ge 1$ be squarefree.
Denote by ${\rm H}_k^*(N)$ the set of all normalized Hecke primitive eigencuspforms
of weight $k$ for the congruence modular group 
$$
\Gamma_0(N)
:=\bigg\{\begin{pmatrix}a & b\\ c & d\end{pmatrix}\in SL_2(\Z) : c\equiv 0\,({\rm mod}\,N)
\bigg\}.
$$ 
Here the normalization is taken to have $\lambda_f(1)=1$ 
in the Fourier series of $f\in {\rm H}_k^*(N)$
at the cusp $\infty$,  
\begin{equation}\label{SFf}
f(z) 
= \sum_{n=1}^\infty \lambda_f(n) n^{(k-1)/2} e^{2\pi inz}
\qquad(\Im m z>0).
\end{equation}
Inherited from the Hecke operators,
the normalized Fourier coefficient $\lambda_f(n)$ satisfies the
following relation
\begin{equation}\label{relationHecke}
\lambda_f(m)\lambda_f(n) 
= \sum_{\substack{d\mid (m, n)\\ (d, N)=1}}
\lambda_f\bigg(\frac{mn}{d^2}\bigg)
\end{equation}
for all integers $m\ge 1$ and $n\ge 1$.
In particular, $\lambda_f(n)$ is multiplicative.
 
Following Deligne \cite{De74},   
for any prime number $p$ 
there are two complex numbers $\alpha_f(p)$ and $\beta_f(p)$ such that
\begin{equation}\label{Deligne1}
\begin{cases}
\alpha_f(p)=\varepsilon_f(p)p^{-1/2}, \; \beta_f(p)=0  & \text{if $p\mid N$}
\\\noalign{\vskip 1mm}
|\alpha_f(p)|=\alpha_f(p)\beta_f(p)=1                  & \text{if $p\nmid N$}
\end{cases}
\end{equation}
and 
\begin{equation}\label{Deligne2}
\lambda_f(p^\nu)
= \frac{\alpha_f(p)^{\nu+1}-\beta_f(p)^{\nu+1}}{\alpha_f(p)-\beta_f(p)}
\end{equation}
for all integers $\nu\ge 1$,
where $\varepsilon_f(p)=\pm 1$.
Hence $\lambda_f(n)$ is real and verifies Deligne's inequality
\begin{equation}\label{Deligne3}
|\lambda_f(n)|
\le d(n)
\end{equation}
for all integers $n\ge 1$, 
where $d(n)$ is the divisor function.
In particular for each prime number $p\nmid N$ there is $\theta_f(p)\in [0,\pi]$ 
such that
\begin{equation}\label{Deligne4}
\lambda_f(p)
= 2\cos\theta_f(p).
\end{equation}
See e.g. \cite{I} for basic analytic facts about modular forms. 

\goodbreak

Positive real moments of Hecke eigenvalues were firstly studied
by Rankin (\cite{Ran83}, \cite{Ran85}).
For $f\in {\rm H}_k^*(N)$ and $r\ge 0$,
consider the sum of the $2r$th power of $|\lambda_f(n)|$:
\begin{equation}\label{defS*}
S_f^*(x; r):=\sum_{n\le x}|\lambda_f(n)|^{2r}.
\end{equation}
The method of Rankin \cite{Ran85} illustrates how to obtain optimally the lower and upper bounds for $S_f^*(x; r)$
if we only know that the associated Dirichlet series
\begin{equation}\label{defDSFr}
F_r(s)
:=\sum_{n\ge 1} |\lambda_f(n)|^{2r}n^{-s}
\quad(\re s>1)
\end{equation}
is  invertible for $\re s\ge 1$ (i.e. holomorphic and nonzero for $\re s\ge 1$) when $r=1, 2$. 
(The invertibility of these two cases are known by Moreno \& Shahidi \cite{MS83}.) 
Rankin's result (\cite{Ran85}, Theorem 1) reads that 
\begin{equation}\label{LUBS*Rankin}
x(\log x)^{\delta_r^{\mp}}
\ll S_f^*(x; r)
\ll x(\log x)^{\delta_r^{\pm}}
\quad
(r\in {\mathcal R}^{\mp})
\end{equation}
for $x\ge x_0(f,r)$,
where 
$$
{\mathcal R}^{-}
:=[0,1]\cup [2,\infty),
\qquad
{\mathcal R}^{+}
:=[1,2],
$$
and
$$
\delta_r^{-}
:=2^{r-1}-1,
\qquad
\delta_r^{+}
:=\frac{2^{r-1}}{5}(2^r+3^{2-r})-1.
$$
The implied constants in (\ref{LUBS*Rankin}) depend on $f$ and $r$.

On the other hand, if the Sato-Tate conjecture holds for newform $f$,
then 
\begin{equation}\label{ST}
S_f^*(x; r)
\sim C_r(f)x(\log x)^{\theta_r}
\quad(x\to\infty),
\end{equation}
where $C_r(f)$ is a positive constant depending on $f, r$ and
$$
\theta_r
:=\frac{4^r\Gamma(r+\frac{1}{2})}{\sqrt{\pi}\Gamma(r+2)}-1.
$$

Very recently,
Tenenbaum \cite{T07} improved
Rankin's exponent 
$
\delta_{1/2}^{+}
=0.0651\cdots
$
to 
$\rho_{1/2}^{+}
=0.1185\cdots
$
(see (\ref{defrho}) below for the definition of $\rho_r^{+}$), 
as an application of his general result 
on the mean values of multiplicative functions and the fact that 
$F_3(s)$ and $F_4(s)$ are invertible for $\re s\ge 1$, 
proven in the excellent work of Kim \& Shahidi \cite{KS02b}.
Although the result (\cite{T07}, Corollary) is stated only for Ramanujan's $\tau$-function,
it is apparent that Tenenbaum's method applies 
to establish the upper bound for $S_f^*(x;r)$ in (\ref{LUBS*}) below.
It should be pointed out that Tenenbaum's approach is different from that of Rankin
and does not give a lower bound for $S_f^*(x;r)$.

The first aim of this paper is to improve the lower and upper bounds
in (\ref{LUBS*Rankin}),
by generalizing Rankin's method to incorporate 
the aforementioned results of Kim \& Shahidi on $F_3(s)$ and $F_4(s)$. 

\begin{theorem}\label{theoremS*}
For any $f\in {\rm H}_k^*(N)$, we have
\begin{equation}\label{LUBS*}
x(\log x)^{\rho_r^{\mp}}
\ll S_f^*(x; r)
\ll x(\log x)^{\rho_r^{\pm}}
\quad
(r\in {\mathscr R}^{\mp})
\end{equation}
for $x\ge x_0(f,r)$,
where 
\begin{equation}\label{defMathscrR}
{\mathscr R}^{-}:=[0,1]\cup [2,3]\cup [4, \infty),
\qquad
{\mathscr R}^{+}:=[1,2]\cup [3,4],
\end{equation}
and
\begin{equation}\label{defrho}
\begin{cases}
\displaystyle 
\rho_r^{-}
:=\frac{3^{r-1}-1}{2},
\\\noalign{\vskip 1mm}
\displaystyle 
\rho_r^{+}
:=\frac{102+7\sqrt{21}}{210}\bigg(\frac{6-\sqrt{21}}{5}\bigg)^r
+\frac{102-7\sqrt{21}}{210}\bigg(\frac{6+\sqrt{21}}{5}\bigg)^r
+\frac{4^r}{35}-1.
\end{cases}
\end{equation}
The implied constants in (\ref{LUBS*}) depend on $f$ and $r$.
\end{theorem}

\smallskip

The upper bound part in (\ref{LUBS*}) 
are essentially due to Tenenbaum \cite{T07},
since his method with a minuscule modification allows to obtain this result.
The lower bound part is new. 
The following table illustrates progress against Rankin's (\ref{LUBS*Rankin})
and the difference from the conjectured values (\ref{ST}).
$$\displaylines{
\vbox{\tabskip = 0pt\offinterlineskip
\halign{
\vrule # & &\hfil$ $ $#$ $ \!\! $ \hfil & \vrule #\cr
\noalign{\hrule}
height 2mm
&&&&&&&&&&&&&&&&&&&&
\cr
& \,\,\, r \, \,\, &
&        0         &
&       0.5        &
& \; 1 \; &
&       1.5        &
& \; 2 \; & 
&       2.5        &
&        3         &
& \,\,\,3.5\, \,\, & 
&        4         &
\cr
height 2mm
&&&&&&&&&&&&&&&&&&&&
\cr
\noalign{\hrule}
height 2mm
&&&&&&&&&&&&&&&&&&&&
\cr
& \delta_r^{-} &
& -0.5         &
& -0.2929      &
& 0            &
& 0.4142       &
& 1            &
& 1.8284       &
& 3            &
& 4.6569       &
& 7            &
\cr
height 2mm
&&&&&&&&&&&&&&&&&&&&
\cr
\noalign{\hrule}
height 2mm
&&&&&&&&&&&&&&&&&&&&
\cr
& \rho_r^{-} &
& -0.3333\,  &
& -0.2113\,  &
& 0          &
& \,0.3660\, &
& 1          &
& \,2.0981\, &
& 4          &
& 7.2945     &
& 13         &
\cr
height 2mm
&&&&&&&&&&&&&&&&&&&&
\cr
\noalign{\hrule}
height 2mm
&&&&&&&&&&&&&&&&&&&&
\cr
& \theta_r  &
& 0         &
& -0.1512   &
& 0         &
& 0.3581    &
& 1         &
& 2.1043    &
& 4         &
& 7.2781    &
& 13        &
\cr
height 2mm
&&&&&&&&&&&&&&&&&&&&
\cr
\noalign{\hrule}
height 2mm
&&&&&&&&&&&&&&&&&&&&
\cr
& \rho_r^{+} &
& 0          &
& -0.1185    &
& 0          &
& 0.3502     &
& 1          &
& 2.1112     &
& 4          &
& 7.2576     &
& 13         &
\cr
height 2mm
&&&&&&&&&&&&&&&&&&&&
\cr
\noalign{\hrule}
height 2mm
&&&&&&&&&&&&&&&&&&&&
\cr
& \delta_r^{+} &
& 0            &
& -0.0652      &
& 0            &
& 0.2899       &
& 1            &
& 2.5266       &
& 5.6667       &
& \,12.0177\,  &
& \,24.7778\,  &
\cr
height 2mm
&&&&&&&&&&&&&&&&&&&&
\cr
\noalign {\hrule}}}
\cr}
$$

\medskip

In order to detect sign changes or cancellations among $\lambda_f(n)$,
it is natural to study summatory function
\begin{equation}\label{SHecke}
S_f(x):=\sum_{n\le x}\lambda_f(n)
\end{equation}
and compare it with (\ref{LUBS*}).
There is a long history on the investigation of the upper estimate for $S_f(x)$.
In 1927, Hecke \cite{He27} showed
$$
S_f(x)\ll_f x^{1/2}
$$
for all $f\in {\rm H}_k^*(N)$ and $x\ge 1$.
Subsequent improvements came with the use of the identity:
$$
\frac{1}{\Gamma(r+1)}
\sum_{n\le x}(x-n)^ra_f(n)
=\frac{1}{(2\pi)^3}\sum_{n\ge 1}
\Big(\frac{x}{n}\Big)^{(k+3)/2}a_f(n)
J_{k+3}\big(4\pi\sqrt{nx}\big),
$$
where $a_f(n):=\lambda_f(n)n^{(k-1)/2}$ and 
$J_{k}(t)$ is the first kind Bessel functions.
Such an identity was first given by Wilton \cite{Wil28}
in which only the case of Ramanujan's $\tau$-function was stated,
and later generalized by Walfisz \cite{Wal33} to other forms.
Let $\vartheta$ be the constant satisfying
$$
|\lambda_f(n)|\ll n^{\vartheta}
\quad(n\ge 1).
$$
Walfisz proved that 
\begin{equation}\label{UBSWalfisz}
S_f(x)\ll_f x^{(1+\vartheta)/3}
\quad(x\ge 1).
\end{equation}
Inserting the values of $\vartheta$ in the historical record into (\ref{UBSWalfisz}) yields
$$
S_f(x)\ll_{f,\varepsilon}
\begin{cases}
x^{11/24+\varepsilon} & \text{Kloosterman \cite{Kl27}}
\\
x^{4/9+\varepsilon}   & \text{Davenport \cite{Da32}, Sali\'e \cite{Sa33}}
\\
x^{5/12+\varepsilon}  & \text{Weil \cite{We48}}
\\
x^{1/3+\varepsilon}   & \text{Deligne \cite{De74}}
\end{cases}
$$
for any $\varepsilon>0$.
Hafner \& Ivi\'c (\cite{HI89}, Theorem 1) removed the factor $x^{\varepsilon}$ of Deligne's result.
On the other hand,  
by combining Walfisz' method with his idea in the study of (\ref{defS*}), 
Rankin \cite{Ran90} showed that  
\begin{equation}\label{UBSRankin}
S_f(x)\ll_{f,\varepsilon} x^{1/3}(\log x)^{\delta_{1/2}^{+}+\varepsilon}
\end{equation}
for any $\varepsilon>0$ and $x\ge 2$.

Here we propose a better bound, by combining Walfisz' method \cite{Wal33} 
and Tenenbaum's approach \cite{T07}.
It is worthy to point out that Tenenbaum's method is 
not only to  improve $\delta_{1/2}^{+}$ to $\rho_{1/2}^{+}$ 
but also remove the $\varepsilon$ in (\ref{UBSRankin}).

\begin{theorem}\label{theoremS}
For $f\in {\rm H}_k^*(N)$, we have
\begin{equation}\label{upperbS}
S_f(x)\ll x^{1/3}(\log x)^{\rho_{1/2}^{+}}
\end{equation}
for $x\ge 2$, where the implied constant depends on $f$.
\end{theorem}

\smallskip

In the opposite direction, Hafner \& Ivi\'c (\cite{HI89}, Theorem 2) proved 
that there is a positive constant $D$ such that
$$
S_f(x)=\Omega_{\pm}\bigg(x^{1/4}\exp\bigg\{\frac{D(\log_2x)^{1/4}}{(\log_3x)^{3/4}}\bigg\}\bigg),
$$
where $\log_r$ denotes the $r$-fold iterated logarithm.

\smallskip

As an application of Theorems \ref{theoremS*} and \ref{theoremS},
we consider the quantities
\begin{equation}\label{LowerNpm}
{\mathscr N}_f^{\pm}(x):=\sum_{\substack{n\le x\\ \lambda_f(n)\gtrless\,0}} 1.
\end{equation}
Very recently Kohnen, Lau \& Shparlinski (\cite{KLS}, Theorem 1) proved
\begin{equation}\label{LowerBN+N-KLS}
{\mathscr N}_f^{\pm}(x)
\gg_f \frac{x}{(\log x)^{17}}
\end{equation}
for $x\ge x_0(f)$.\footnote{It is worthy to indicate that they gave explicit values 
for the implied constant in $\ll$ and $x_0(f)$.}

Here we propose a better bound.

\begin{corollary}\label{LBNpmx}
For any $f\in {\rm H}_k^*(N)$, we have
$$
{\mathscr N}_f^{\pm}(x)
\gg \frac{x}{(\log x)^{1-1/\sqrt{3}}}
$$
for $x\ge x_0(f)$, where the implied constant depends on $f$.
If we assume Sato-Tate's conjecture, 
the exponent $1-1/\sqrt{3}\approx 0.422$ can be improved to $2-16/(3\pi)\approx 0.302$.
\end{corollary}

In a joint paper with Lau \cite{LW08}, 
we shall remove the logarithmic factor by a completely different method.

\medskip

{\bf Acknowledgment}.
The author would like to thank Winfried Kohnen for the preprint \cite{KLS}
and Yuk Kam Lau for his many suggestions that improved the writting of this paper.

\vskip 10mm

\section{Method of Rankin}

\smallskip

Let $k\ge 2$ be an even integer, $N\ge 1$ be squarefree, $f\in {\rm H}_k^*(N)$ and $r>0$.
Following Rankin's idea \cite{Ran85}, 
we shall find two optimal multiplicative functions $\lambda_{f,r}^{\pm}(n)$ such that
\begin{equation}\label{LUBlambdafr}
\lambda_{f,r}^{\mp}(p^\nu)
\le |\lambda_f(p^\nu)|^{2r}
\le \lambda_{f,r}^{\pm}(p^\nu)
\quad(r\in {\mathscr R}^{\mp})
\end{equation}
for all primes $p$ and integers $\nu\ge 1$,
and furthermore, their associated Dirichlet series $\Lambda_{f,r}^{\pm}(s)$
(see (\ref{defDSLambdafr}) below)
in the half-plane $\re s\ge 1$ is controlled by $F_j(s)$ for $j=1, \dots, 4$.
Then we can apply Tauberian theorems to obtain 
the asymptotic behaviour of the summatory functions of $\lambda_{f,r}^{\pm}(n)$.

\subsection{Construction of $\lambda_{f,r}^{+}(n)$}

\smallskip

For $\boldsymbol{a}:=(a_1, \dots, a_4)\in \R^4$ and $r>0$,
consider the function
\begin{equation}\label{defhr}
h_r(t; \boldsymbol{a}):=t^r-a_1t-a_2t^2-a_3t^3-a_4t^4
\qquad(0\le t\le 1)
\end{equation}
and let
\begin{equation}\label{kappaeta}
\kappa_{-}:=\textstyle\frac{1}{4},
\qquad
\eta_{-}:=\frac{3}{4},
\qquad
\kappa_{+}:=\frac{6-\sqrt{21}}{20},
\qquad
\eta_{+}:=\frac{6+\sqrt{21}}{20}.
\end{equation}
In Subsection \ref{2.3}, 
we shall explain the reason behind this choice.

\medskip

\begin{lemma}\label{a-}
If the function $h_r(t; \boldsymbol{a})$ defined by (\ref{defhr}) satisfies
$$h_r'(\kappa_{-}; \boldsymbol{a})
=h_r'(\eta_{-}; \boldsymbol{a})
=h_r(\kappa_{-}; \boldsymbol{a})
=h_r(\eta_{-}; \boldsymbol{a})
=0,
$$
then 
\begin{equation}\label{defaj-}
a_j
=a_j^{-}
:=\frac{P_j^{-}(\kappa_{-},\eta_{-})-P_j^{-}(\eta_{-},\kappa_{-})}
{(\kappa_{-}-\eta_{-})^3}
\end{equation} 
for $1\le j\le 4$, where
\begin{align*}
P_1^{-}(\kappa,\eta)
& :=\{(4-r)\kappa+(r-2)\eta\}\kappa^{r-1}\eta^2,
\\
P_2^{-}(\kappa,\eta)
& :=\{(2r-8)\kappa^2+(1-r)\kappa\eta+(1-r)\eta^2\}
\kappa^{r-2}\eta, 
\\
P_3^{-}(\kappa,\eta)
& :=\{(4-r)\kappa^2+(4-r)\kappa\eta+2(r-1)\eta^2\}
\kappa^{r-2},
\\
P_4^{-}(\kappa,\eta)
& :=\{(r-3)\kappa+(1-r)\eta\}\kappa^{r-2}.
\end{align*}
\end{lemma}

\proof
This can be done by routine calculation.
\hfill
$\square$

\medskip

\begin{lemma}\label{a+}
If the function $h_r(t; \boldsymbol{a})$ defined by (\ref{defhr}) is such that
$$
\begin{cases}
h_r'(\kappa_{+}; \boldsymbol{a})
=h_r'(\eta_{+}; \boldsymbol{a})=0,
\\\noalign{\vskip 1mm}
h_r(\kappa_{+}; \boldsymbol{a})
=h_r(\eta_{+}; \boldsymbol{a})
=h_r(1; \boldsymbol{a}),
\end{cases}
$$
then 
\begin{equation}\label{defaj+}
a_j
=a_j^{+}
:=\frac{P_j^{+}(\kappa_{+},\eta_{+})-P_j^{+}(\eta_{+},\kappa_{+})}
{(\kappa_{+}-1)^2(\eta_{+}-1)^2(\kappa_{+}-\eta_{+})^3}
\end{equation}
for $1\le j\le 4$, where
\begin{align*}
P_1^{+}(\kappa,\eta)
& :=r\kappa^{r-1}\eta(\kappa-1)(\eta-\kappa)(\kappa\eta+2\kappa+\eta)(\eta-1)^2
\\
& \quad
+2(\kappa^r-1)\kappa\eta(\eta-1)^2(2\kappa\eta+4\kappa-\eta^2-2\eta-3),
\\\noalign{\vskip 1mm}
P_2^{+}(\kappa,\eta)
& :=r\kappa^{r-1}(\kappa-1)(\kappa-\eta)(\eta-1)^2(2\kappa\eta+\kappa+\eta^2+2\eta)
\\
& \quad
+(\eta^r-1)(\kappa-1)^2
(8\kappa\eta^2+4\eta^2-\eta\kappa^2-2\kappa\eta-3\eta-\kappa^3-2\kappa^2-3\kappa), 
\\\noalign{\vskip 1mm}
P_3^{+}(\kappa,\eta)
& :=r\kappa^{r-1}(\kappa-1)(\kappa+2\eta+1)(\eta-\kappa)(\eta-1)^2
\\
& \quad
+2(\kappa^r-1)(2\kappa^2+2\kappa\eta-\eta^2-2\eta-1)(\eta-1)^2,
\\\noalign{\vskip 1mm}
P_4^{+}(\kappa,\eta)
& :=r\kappa^{r-1}(\kappa-1)(\kappa-\eta)(\eta-1)^2
+(\eta^r-1)(\kappa-1)^2(3\eta-\kappa-2).
\end{align*}
\end{lemma}

\proof
This is done by routine calculation as well.
\hfill
$\square$

\medskip

\begin{lemma}\label{LUBoundhr}
Let $\boldsymbol{a}^{\pm}:=(a_1^{\pm}, \dots, a_4^{\pm})$,
where each $a_i^{\pm}$ is given by the value in Lemmas \ref{a-}-\ref{a+}, respectively.
Then for $0\le t\le 1$ we have
$$
h_r(t; \boldsymbol{a}^-)\gtrless 0
\quad\hbox{and}\quad
h_r(t; \boldsymbol{a}^+)\lessgtr h_r(1; \boldsymbol{a}^+)
\quad\hbox{for}\quad
r\in {\mathscr R}^{\mp}.
$$
\end{lemma}

\proof
We have
$$
h_r^{(4)}(t; \boldsymbol{a}^{-})
=r(r-1)(r-2)(r-3)t^{r-4}-24a_4^{-},
$$
so $h_r^{(4)}(t; \boldsymbol{a}^{-})$ has at most one zero for $t>0$
and $h_r^{(i)}(t; \boldsymbol{a}^{-})$ has at most $5-i$ zeros for $t>0$ 
($i=3, 2, 1, 0$).
Since 
$h_r(\kappa_{-}; \boldsymbol{a}^{-})
=h_r(\eta_{-}; \boldsymbol{a}^{-})
=h_r(0; \boldsymbol{a}^{-})$,
it follows that 
$h_r'(\xi_{-}; \boldsymbol{a}^{-})
=h_r'(\xi_{-}'; \boldsymbol{a}^{-})
=0$
for some $\xi_{-}\in (0, \kappa_{-})$
and $\xi_{-}'\in (\kappa_{-}, \eta_{-})$.
Therefore $\xi_{-}$, $\kappa_{-}$, $\xi_{-}'$ and $\eta_{-}$
are the only zeros of $h_r'(t; \boldsymbol{a}^{-})$ in $(0, 1)$.

Now 
$$
h_r''(\kappa_{-}; \boldsymbol{a}^{-})
=8\cdot 4^{-r}
(2r^2-2r+3+2r3^{r-2}-11\cdot 3^{r-2})
$$
and
$$
h_r''(\eta_{-}; \boldsymbol{a}^{-})
=8\cdot 4^{-r}
(2r^2-6r-3-2r3^r+43\cdot 3^{r-2}).
$$
From these, it is easy to verify that
$$
h_r''(\kappa_{-}; \boldsymbol{a}^{-}),
\;
h_r''(\eta_{-}; \boldsymbol{a}^{-})
\begin{cases}
\gtrless 0 & \text{if} \;\, r\in {{\mathscr R}^{\mp}}^{^{\hskip -4,3mm\circ}}\;\;,
\\\noalign{\vskip 1mm}
=0 & \text{if} \;\, r=1, 2, 3, 4,
\end{cases}
$$
where ${{\mathscr R}^{\mp}}^{^{\hskip -4,3mm\circ}}\;\;$ denotes the interior of ${\mathscr R}^{\mp}$.
Hence $h_r(t; \boldsymbol{a}^{-})$ takes its mimimum 
(maximum, respectively) values in $[0, 1]$
at $0$, $\kappa_{-}$, $\eta_{-}$
when $r\in {{\mathscr R}^{-}}^{^{\hskip -4,3mm\circ}}\;\;$ 
($r\in {{\mathscr R}^{+}}^{^{\hskip -4,3mm\circ}}\;\;$, respectively).
Moreover, $h_r(t; \boldsymbol{a}^{-})$ has local maxima
(minima, respectively) at $\xi_{-}$, $\xi_{-}'$
when $r\in {{\mathscr R}^{-}}^{^{\hskip -4,3mm\circ}}\;\;$ 
($r\in {{\mathscr R}^{+}}^{^{\hskip -4,3mm\circ}}\;\;$, respectively).
This proves the assertion about $h_r(t; \boldsymbol{a}^{-})$.

Similarly we can prove the corresponding result on 
$h_r(t; \boldsymbol{a}^{+})$.
\hfill
$\square$

\medskip

Now we define the multiplicative function $\lambda_{f,r}^{\pm}(n)$ by
\begin{equation}\label{deflambdarmp}
\lambda_{f,r}^{\mp}(p^\nu)
:= \begin{cases}
\sum_{0\le j\le 4} 2^{2(r-j)}a_j^{\mp}\lambda_f(p)^{2j} & \text{if $\;\nu=1$ and $r>0$},
\\\noalign{\vskip 1mm}
0                       & \text{if $\;\nu\ge 2$ and $r\in {\mathscr R}^{\mp}$},
\\\noalign{\vskip 1mm}
|\lambda_f(p^\nu)|^{2r} & \text{if $\;\nu\ge 2$ and $r\in {\mathscr R}^{\pm}$},
\end{cases}
\end{equation}
where
\begin{equation}\label{defa0mp}
a_0^{-}:=0
\quad{\rm and}\quad
a_0^{+}:=1-a_1^{+}-a_2^{+}-a_3^{+}-a_4^{+}.
\end{equation}
In view of (\ref{Deligne4}), 
we can apply Lemma \ref{LUBoundhr} with $t=|\cos\theta_f(p)|$ to deduce that
the inequality (\ref{LUBlambdafr}) hold for all primes $p$ and integers $\nu\ge 1$.
Thanking to the multiplicativity, 
these inequalities also hold for all integers $n\ge 1$.

\subsection{Dirichlet series associated to $\lambda_{f,r}^{\pm}(n)$}
For $f\in {\rm H}_k^*(N)$, $r>0$ and $\re s>1$, we define
\begin{equation}\label{defDSLambdafr}
\Lambda_{f,r}^{\pm}(s)
:=\sum_{n\ge 1} \lambda_{f,r}^{\pm}(n)n^{-s}.
\end{equation}
Next we shall study their analytic properties in the half-plane $\re s\ge 1$
by using the higher order symmetric power $L$-functions $L(s, {\rm sym}^mf)$
associated to $f\in {\rm H}_k^*(N)$,
due to 
Gelbart \& Jacquet \cite{GJ78} for $m=2$,
Kim \& Shahidi (\cite{KS02a}, \cite{KS02b}) for $m=3, 4, 5, 6, 7, 8$.
Here the symmetric $m$th power associated to $f$ is defined as
$$ 
L(s, {\rm sym}^mf)
:= \prod_p\prod_{0\le j\le m}
\big(1-\alpha_f(p)^{m-j}\beta_f(p)^{j}p^{-s}\big)^{-1} 
$$
for $\re s>1$, where $\alpha_f(p)$ and $\beta_f(p)$ are given by (\ref{Deligne1}) and (\ref{Deligne2}).
According to the literature mentioned above,
it is known that
the function $L(s, {\rm sym}^mf)$ for $m=2, 3, \dots, 8$
is invertible for $\re s\ge 1$.

\smallskip

We start to study $F_{1}(s)$, $F_{2}(s)$, $F_{3}(s)$ and $F_{4}(s)$.

\begin{lemma}\label{F1234s}
Let $k\ge 2$ be an even integer, $N\ge 1$ be squarefree and $f\in {\rm H}_k^*(N)$.
For $j=1, 2, 3, 4$ and $\re s>1$, we have
\begin{equation}\label{Expression}
F_j(s)=\zeta(s)^{m_j}G_j(s)H_j(s),
\end{equation}
where
\begin{equation}\label{defmj}
m_1:=1,
\qquad
m_2:=2,
\qquad
m_3:=5,
\qquad
m_4:=14,
\end{equation}
and
\begin{align*}
G_1(s)
& := L(s, {\rm sym}^2f),
\\ 
G_2(s)
& := L(s, {\rm sym}^2f)^3L(s, {\rm sym}^4f),
\\
G_3(s)
& := L(s, {\rm sym}^2f)^9L(s, {\rm sym}^4f)^5L(s, {\rm sym}^6f),
\\
G_4(s)
& := L(s, {\rm sym}^2f)^{34}L(s, {\rm sym}^4f)^{20}
L(s, {\rm sym}^6f)^7L(s, {\rm sym}^8f)
\end{align*}
are invertible for $\re s\ge 1$.
Here the function $H_j(s)$ admits a Dirichlet series convergent absolutely in $\re s>\frac{1}{2}$ and $H_j(s)\not=0$ for $\re s=1$.
\end{lemma}

\proof
Write $x$ for the trace of a local factor of $L(s,f)$ (i.e. $\alpha_f(p)+\beta_f(p)$), and 
denote by  $T_n(x)$ the polynomial 
for the trace of its symmetric $n$th power.
Then 
\begin{align*}
T_2
& =x^2-1,
\\
T_4
& =x^4-3x^2+1,
\\
T_6
& =x^6-5x^4+6x^2-1,
\\
T_8
& =x^8-7x^6+15x^4-10x^2+1,
\end{align*}
from which we deduce
\begin{align*}
x^2
& =1+T_2,
\\
x^4
& =2+3T_2+T_4,
\\
x^6
& =5+9T_2+5T_4+T_6,
\\
x^8
& =14+34T_2+20T_4+7T_6+T_8.
\end{align*}
This implies (\ref{Expression}).
By using results on $L(s, {\rm sym}^mf)$ mentioned above,
$G_j(s)$ is invertible for $\re s\ge 1$.
This completes the proof.
\hfill
$\square$

\smallskip

\begin{lemma}\label{Lambdarpms}
Let $k\ge 2$ be an even integer, $N\ge 1$ be squarefree and $f\in {\rm H}_k^*(N)$.
For $r>0$ and $\re s>1$, we have
\begin{equation}\label{ExpLambda}
\Lambda_{f,r}^{\pm}(s)
= \zeta(s)^{\rho_r^{\pm}+1}H_{f,r}^{\pm}(s),
\end{equation}
where
\begin{equation}\label{defrhorpm}
\rho_r^{\pm}:=2^{2r-8}(2^{8}a_0^{\pm}+2^{6}a_1^{\pm}+2^{4}\cdot 2a_2^{\pm}
+2^{2}\cdot 5a_3^{\pm}+14a_4^{\pm})-1
\end{equation}
and $H_{f,r}^{\pm}(s)$ is invertible for $\re s\ge 1$.
\end{lemma}

\proof
By definition (\ref{deflambdarmp}),
for $\re s>1$ we can write
\begin{align*}
\Lambda_{f,r}^{-}(s)
& = \prod_p 
\Big(1+\sum_{0\le j\le 4} 2^{2(r-j)}a_j^{-}\lambda_f(p)^{2j}p^{-s}\Big)
\\
& = \prod_{0\le j\le 4} F_j(s)^{2^{2(r-j)}a_j^{-}}
H_r^{-}(s)
\end{align*}
for $r\in {\mathscr R}^{-}$, and
\begin{align*}
\Lambda_{f,r}^{-}(s)
& = \prod_p 
\Big(1+\sum_{0\le j\le 4} 2^{2(r-j)}a_j^{-}\lambda_f(p)^{2j}p^{-s}
+\sum_{\nu\ge 2} |\lambda_f(p^\nu)|^{2r}p^{-\nu s}\Big)
\\
& = \prod_{0\le j\le 4} F_j(s)^{2^{2(r-j)}a_j^{-}}
H_r^{-}(s)
\end{align*}
for $r\in {\mathscr R}^{+}$,
where $F_0(s)=\zeta(s)$ is the Riemann zeta-function and
$H_r^{-}(s)$ is a Dirichlet series absolutely convergent for $\re s>\frac{1}{2}$
such that $H_r^{-}(s)\not=0$ for $\re s=1$.
Now the desired result with the sign `$-$' follows from Lemma \ref{F1234s}.
The other part can be treated in the same way.
\hfill
$\square$

\subsection{Optimalisation of $\lambda_{f,r}^{\pm}(p)$
and choice of $\kappa_{\pm}, \eta_{\pm}$}\label{2.3}
If we regard $\kappa_{\pm}, \eta_{\pm}$ as parameters,
the $\rho_r^{\pm}$ given by (\ref{defrhorpm}) are functions of these parameters.
We choose $(\kappa_{\pm}, \eta_{\pm})$ in $(0,1)^2$ optimally,
which can be done by using formal calculation via Maple.
Their values are given by (\ref{kappaeta}).

\vskip 10mm

\section{Proof of Theorem \ref{theoremS*}}

\smallskip

In view of Lemma \ref{Lambdarpms} and the classical fact on $\zeta(s)$,
we can write
\begin{equation}\label{neighbourhood}
\Lambda_{f,r}^{\pm}(s)
=\frac{H_{f,r}^{\pm}(1)}{(s-1)^{\rho_r^{\pm}+1}}
+g_{f,r}^{\pm}(s)
\end{equation}
in some neighbourhood of $s=1$ with $\re s>1$,
where $H_{f,r}^{\pm}(1)\not=1$ and $g_{f,r}^{\pm}(s)$ is holomorphic at $s=1$.
Since $\lambda_{f,r}^{\pm}(n)\ge 0$,
we can apply Delange's tauberian theorem \cite{Delange54} to write
\begin{equation}\label{asymp}
\sum_{n\le x} \lambda_{f,r}^{\pm}(n)
\sim H_{f,r}^{\pm}(1)x(\log x)^{\rho_r^{\pm}}
\quad(x\to\infty).
\end{equation}
Now Theorem \ref{theoremS*} follows from (\ref{LUBlambdafr}) and (\ref{asymp}).

\vskip 10mm

\section{Proof of Theorem \ref{theoremS}}

\smallskip

By (\ref{neighbourhood}), it follows that
$$
\prod_p\bigg(1+\sum_{\nu\ge 1} \frac{\lambda_{f,r}^{\pm}(p^\nu)}{p^{\nu\sigma}}\bigg)
=\frac{H_{f,r}^{\pm}(1)}{(\sigma-1)^{\rho_r^{\pm}+1}}
+g_{f,r}^{\pm}(\sigma)
$$
for $\sigma>1$.
From this, (\ref{deflambdarmp}), (\ref{defa0mp}) and Deligne's inequality, we deduce that
$$
\sum_p  \frac{\lambda_{f,r}^{\pm}(p)}{p^{\sigma}}
=(\rho_r^{\pm}+1)\log(\sigma-1)^{-1}+C_{f,r}^{\pm}+o(1)
\quad(\sigma\to 1+),
$$
where $C_{f,r}^{\pm}$ is some constant.

On the other hand, the prime number theorem implies, by a partial integration, that
$$
\sum_p p^{-\sigma}=\log(\sigma-1)^{-1}+C+o(1)
\quad(\sigma\to 1+),
$$
where $C$ is an absolute constant.
Thus the preceding relation can be written as
\begin{equation}\label{TeEx78}
\sum_p  \frac{\lambda_{f,r}^{\pm}(p)-(\rho_r^{\pm}+1)}{p^{\sigma}}
=C_{f,r}^{\pm}+(\rho_r^{\pm}+1)C+o(1)
\quad(\sigma\to 1+).
\end{equation}
According to Exercise II.7.8 of \cite{T96},
the formula (\ref{TeEx78}) implies
$$
\sum_p \frac{\lambda_{f,r}^{\pm}(p)-(\rho_r^{\pm}+1)}{p}
=C_{f,r}^{\pm}+(\rho_r^{\pm}+1)C.
$$
Hence 
$$
\sum_{p\le x}  \frac{\lambda_{f,r}^{\pm}(p)}{p}
=(\rho_r^{\pm}+1)\log_2x+C_{f,r}^{\pm}+(\rho_r^{\pm}+1)C+o(1)
\quad(x\to\infty).
$$

Now we apply a well known result of Shiu \cite{Sh80} and (\ref{LUBlambdafr}) to write
\begin{equation}\label{Shiu}
\begin{aligned}
\sum_{x\le n\le x+z}|\lambda_f(n)|^{2r}
& \ll \frac{z}{\log x}\exp\bigg(\sum_{p\le x} \frac{|\lambda_f(p)|^{2r}}{p}\bigg)
\\
& \ll \frac{z}{\log x}\exp\bigg(\sum_{p\le x} \frac{\lambda_{f,r}^{+}(p)}{p}\bigg)
\\
& \ll z(\log x)^{\rho_{r}^{+}}
\end{aligned}
\end{equation}
for $r\in {\mathscr R}^{-}$, any $\varepsilon>0$, $x\ge x_0(\varepsilon)$ and $x^{1/4}\le z\le x$.
Using this with $r=\frac{1}{2}$ in (9) of \cite{Ran90},
the first term on the right-hand side of (10) of \cite{Ran90} is replaced by
$x^{1/2}z^{-1/2}(\log x)^{\rho_{1/2}^{+}}$.
Applying (\ref{Shiu}) with $r=\frac{1}{2}$ again to the second term 
on the right-hand side of (10) of \cite{Ran90}, it follows that
$$
S_f(x)
\ll x^{1/2}z^{-1/2}(\log x)^{\rho_{1/2}^{+}}
+z(\log x)^{\rho_{1/2}^{+}}.
$$
Taking $z=x^{1/3}$,
we obtain the required result when the level is $N=1$.
The general case can be treated similarly as indicated in \cite{Ran90}.
\hfill
$\square$

\vskip 10mm

\section{Proof of Corollary \ref{LBNpmx}}

\smallskip

By comparing (\ref{upperbS}) and the lower bound part in (\ref{LUBS*}) with $r=\frac{1}{2}$,
it is easy to deduce that 
$$
\sum_{\substack{n\le x\\ \lambda_f(n)\gtrless\,0}} |\lambda_f(n)|
\gg_f x(\log x)^{\rho_{1/2}^{-}}
$$
for $x\ge x_0(f)$.
Since $\rho_{1/2}^{-}=-(1-1/\sqrt{3})/2$ and $\rho_{1}^{+}=0$,
a simple application of the Cauchy-Schwarz inequality yields the following result.

The second assertion can be obtained by noticing that 
$\theta_{1/2}=8/(3\pi)-1$.
\hfill
$\square$

\end{document}